\documentclass[a4paper,11 pt]{amsart}

\RequirePackage{amsmath, amssymb, amsthm, amsfonts}
\usepackage{fullpage}
\newtheorem{theo}{Theorem}[section]

\newtheorem{prop}[theo]{Proposition}

\newtheorem{defin}[theo]{Definition}

\theoremstyle{definition}

\newtheorem{rem}[theo]{Remark}

\newcommand{\PP}{\mathbb{P}}

\newcommand{\oo}{\mathcal{O}}

\newcommand{\Aut}{{\mathrm{Aut}}}

\newcommand{\ord}{{\mathrm{ord}}}

\newcommand{\ZZ}{\mathbb{Z}}

\title{On asymptotic bounds for the number of\\ irreducible components of the moduli space of\\ surfaces of general type II}
\author{Michael L\"onne and Matteo Penegini}
\address{Michael L\"onne\\  
Mathematik VIII, Universit\"at Bayreuth,
NWII, Universit\"atsstrasse 30, D-95447 Bayreuth,  Germany}
\email{michael.loenne@uni-bayreuth.de}
\address{Matteo Penegini\\
Dipartimento di Matematica \emph{``Federigo Enriques''}, Universit\`{a} degli Studi di Milano, Via Saldini 50, I-20133 Milano, Italy} \email{matteo.penegini@unimi.it}

\begin{document}


\maketitle


\begin{abstract}
In this paper we investigate the asymptotic growth of the number of irreducible and connected components of the moduli space of
surfaces of general type corresponding to certain families of
surfaces isogenous to a higher product with group $(\ZZ/2\ZZ)^k$. We obtain a significantly higher growth than the one in our previous paper  \cite{LP14}.
\end{abstract}

\section{Introduction}

It is well known (see \cite{Gie77}) that once two positive integers $x,y$ are fixed there exists a quasiprojective coarse moduli space $\mathcal{M}_{y,x}$ of canonical models of surfaces of general type with $x=\chi(S)$ and $y=K^2_S$. The number $\iota(x,y)$, resp.\ $\gamma(x,y)$, of irreducible, resp.\ connected, components of  $\mathcal{M}_{y,x}$ is bounded from above by a function of $y$.
In fact, Catanese proved that the number $\iota^0(y,x)$ of components containing regular surfaces, i.e., \ $q(S)=0$, has
an exponential upper bound in $K^2$. More precisely \cite[p.592]{Cat92} gives the following inequality
\[
\iota^0(x,y) \leq y^{77y^2}. 
\]
This result is not known to be sharp and in recent papers \cite{M97,Ch96,GP14,LP14} inequalities are proved which tell how close one can get to this bound from below. In particular, in the last two papers the authors considered families of surfaces isogenous to a product in order to construct many irreducible components of the moduli space of surfaces of general type. The reason why one works with these surfaces, is the fact that the number of families of these surfaces can be easily computed using group theoretical and combinatorial methods.

In our previous work \cite{LP14} we constructed many such families with many different $2-$groups. There, we esploited the fact that the number of $2-$groups with given order grows very fast in function of the order. In this paper we obtain a significantly better lower bound for $\iota^0(x,y)$ using only the groups $(\ZZ/2\ZZ)^k$ and again some properties of the moduli space of surfaces isogenous to a product. Our main result is the following theorem.

\begin{theo}\label{thm main}\sloppy Let $h$ be number of connected components of the moduli space of surfaces of general type containing regular surfaces isogenous to a product of curves, admitting $(\ZZ/2\ZZ)^k$ as group and ramification structure of type $(2^{k(k-1)/2},2^{k^2-k-1})$. Then for $k \rightarrow \infty$ we have  
\[ h \geq 2^{c \sqrt[2+\nu]{x_k}}.
\]
with $\nu$ a positive real number. In particular, we obtain sequences $y_k$ with
\[\iota^0(x_k,y_k) \geq Cy_k^{\sqrt{y_k}}.
\]
\end{theo}
\noindent
Let us explain now the way in which this paper is organized.

In the next section \emph{Preliminaries}  we recall the definition and the properties of surfaces isogenous to a higher product and the its associated group theoretical data. Moreover, we recall a result of Bauer--Catanese \cite{BC} which allows us to count the number of connected components of the moduli space of surfaces isogenous to a product with given group and type of ramification structure.

In the last section we give the proof of the Theorem \ref{thm main}.

\textbf{Acknowledgement} The first author was  supported by  the ERC 2013 Advanced Research Grant - 340258 - TADMICAMT,
at the Universitaet Bayreuth. The second author acknowledges Riemann Fellowship Program of the Leibniz Universit\"at Hannover. 

\bigskip

\textbf{Notation and conventions.} We work over the field
$\mathbb{C}$ of complex numbers. By \emph{surface} we mean a projective, non-singular surface $S$. For such a surface $\omega_S=\oo_S(K_S)$ denotes the canonical
bundle, $p_g(S)=h^0(S, \, \omega_S)$ is the \emph{geometric genus},
$q(S)=h^1(S, \, \omega_S)$ is the \emph{irregularity},
$\chi(\mathcal{O}_S)=\chi(S)=1-q(S)+p_g(S)$ is the \emph{Euler-Poincar\'e
characteristic} and $e(S)$ is the \emph{topological Euler number} of $S$. 
\section{Preliminaries}
\begin{defin}\label{def.isogenous} A surface $S$ is said to be \emph{isogenous to a higher product of curves}\index{Surface isogenous to a higher product of curves} if and only if, $S$ is a
quotient $(C_1 \times C_2)/G$, where $C_1$ and $C_2$ are curves of
genus at least two, and $G$ is a finite group acting freely on
$C_1 \times C_2$.
\end{defin}
Using the same notation as in Definition \ref{def.isogenous}, let
$S$ be a surface isogenous to a product, and $G^{\circ}:=G
\cap(Aut(C_1) \times Aut(C_2))$. Then $G^{\circ}$ acts on the two
factors $C_1$, $C_2$ and diagonally on the product $C_1 \times
C_2$. If $G^{\circ}$ acts faithfully on both curves, we say that
$S= (C_1 \times C_2)/G$ is a \emph{minimal
realization}. In \cite{cat00}
it is also proven that any
surface isogenous to a product
admits a unique minimal realization. 

\medskip

{\bf Assumptions.} In the following we always assume:
\begin{enumerate}
\item Any surface $S$ isogenous to a product is given by its unique minimal realization;
\item $G^{\circ}=G$, this case is also known as \emph{unmixed type}, see \cite{cat00}.
\end{enumerate}
Under these assumption we have. \nopagebreak
\begin{prop}~\cite{cat00}\label{isoinv}
Let $S=(C_1 \times C_2)/G$ be a surface isogenous to a higher product of curves, then $S$ is a minimal surface of general type with the following invariants:
\begin{equation}\label{eq.chi.isot.fib}
\chi(S)=\frac{(g(C_1)-1)(g(C_2)-1)}{|G|},
\quad
e(S)=4 \chi(S),
\quad
K^2_S=8 \chi(S).
\end{equation}
The irregularity of these surfaces is computed by
\begin{equation}\label{eq_irregIsoToProd}
 q(S)=g(C_1/G)+g(C_2/G).
\end{equation}
\end{prop}
Among the nice features of surfaces isogenous to a product, one is that their deformation class can be obtained in a purely algebraic way. Let us briefly recall this in the particular case when $S$ is regular, i.e., $q(S)=0$, $C_i/G \cong \PP^1$. 
\begin{defin}\label{defn.sphergen}
Let $G$ be a finite group and $r \in \mathbb{N}$ with $r \geq 2$.
\begin{itemize}
\item An $r-$tuple $T=(v_1,\ldots,v_r)$ of elements of $G$ is called a
\emph{spherical system of generators} of $G$ if $ \langle
v_1,\ldots,v_r \rangle=G$ and $v_1 \cdot \ldots \cdot v_r=1$.

\item We say that $T$ has an \emph{unordered type} $\tau:=(m_1, \dots ,m_r)$ if the
orders of $(v_1,\dots,v_r)$ are $(m_1, \dots ,m_r)$ up to a
permutation, namely, if there is a permutation $\pi \in \mathfrak{S}_r$ such
that
\[
    \ord(v_1) = m_{\pi(1)},\dots,\ord(v_r)= m_{\pi(r)}.
\]

\item Moreover, two spherical systems $T_1=(v_{1,1}, \dots , v_{1,r_1})$ and $T_2=(v_{2,1}, \dots , w_{2,r_2})$
are said to be \emph{disjoint}, if:
\begin{equation}\label{eq.sigmasetcond} \Sigma(T_1)
\bigcap \Sigma(T_2)= \{ 1 \},
\end{equation}
where
\[ \Sigma(T_i):= \bigcup_{g \in G} \bigcup^{\infty}_{j=0} \bigcup^{r_i}_{k=1} g \cdot v^j_{i,k} \cdot
g^{-1}.
\]
\end{itemize}
\end{defin}
We shall also use the shorthand, for example $(2^4,3^2)$, to indicate
the tuple $(2,2,2,2,3,3)$. 
\begin{defin}\label{def.rami.structure} Let $2<r_i \in \mathbb{N}$ for $i=1,2$ and $\tau_i=(m_{i,1}, \dots ,m_{i,r_i})$ be two sequences of natural numbers such
that $m_{k,i} \geq 2$.  A \emph{(spherical-) ramification structure} of type $(\tau_1,\tau_2)$ and size
$(r_1,r_2)$ for a finite group $G$, is a
pair $(T_1,T_2)$
of disjoint spherical systems of generators of $G$, whose types are
$\tau_i$, such that:
\begin{equation}\label{eq.Rim.Hur.Condition}
\mathbb{Z} \ni \frac{|G|
(-2+\sum^{r_i}_{l=1}(1-\frac{1}{m_{i,l}}))}{2}+1 \geq 2,\qquad \text{for } i=1,2.
\end{equation}
\end{defin}
\begin{rem}\label{minimal}
Following e.g., the discussion in \cite[Section 2]{LP14} we obtain that the datum of the deformation class of a regular surface
$S$ isogenous to a higher product of curves of unmixed type together with its minimal realization $S=(C_1 \times
C_2)/G$ is determined by the datum of a finite
group $G$ together with two
disjoint spherical systems of generators $T_1$ and $T_2$ (for more details see also \cite{BCG06}).
\end{rem}
\begin{rem} Recall that from Riemann Existence Theorem a finite group $G$ acts as a
group of automorphisms of some curve $C$ of
genus $g$ such that $C/G \cong \PP^1$ if and only if there exist integers $m_r
\geq m_{r-1} \geq \dots \geq m_1\geq 2$ such that $G$ has a spherical system of generators of type $(
m_1,\dots,m_r)$ and the following
Riemann-Hurwitz relation holds:
\begin{equation}\label{eq.RiemHurw} 2g-2=| G | (-2 +
\sum^r_{i=1}(1-\frac{1}{m_i})).
\end{equation}
\end{rem}
\begin{rem} Note that a group $G$ and a ramification structure determine the main numerical
invariants of the surface $S$. Indeed, by \eqref{eq.chi.isot.fib} and~\eqref{eq.RiemHurw} we obtain:
\begin{equation}\label{eq.pginfty}
4\chi(S)=|G|\cdot\left({-2+\sum^{r_1}_{k=1}(1-\frac{1}{m_{1,k}})}\right)
\cdot\left({-2+\sum^{r_2}_{k=1}(1-\frac{1}{m_{2,k}})}\right)=: 4\chi (|G|,(\tau_1,\tau_2)).
\end{equation}
\end{rem}
Let $S$ be a surface isogenous to a product of unmixed type with group $G$ and a pair of
two disjoint spherical systems of generators of types
$(\tau_1,\tau_2)$. By~$\eqref{eq.pginfty}$ we have
$\chi(S)=\chi(G,(\tau_1,\tau_2))$, and consequentially,
by~\eqref{eq.chi.isot.fib}, $K^2_S=K^2(G,(\tau_1,\tau_2))=8\chi(S)$.

Let us fix a group $G$ and a pair of unmixed ramification types
$(\tau_1,\tau_2)$, and denote by $\mathcal{M}_{(G,(\tau_1,\tau_2))}$
the moduli space of isomorphism classes of surfaces isogenous to a product
admitting these data, by
\cite{cat00,cat04} the space $\mathcal{M}_{(G,(\tau_1,\tau_2))}$ consists
of a finite number of connected components. Indeed, there is a group
theoretical procedure to count these components. In case $G$ is abelian it is described in
\cite{BC}.

\begin{theo}\cite[Theorem 1.3]{BC} \label{theo_count}.
Let $S$ be a surface isogenous to a higher product
of unmixed type and with $q=0$.
Then to $S$ we attach its finite group $G$ (up to isomorphism) and
the equivalence classes of an unordered pair of disjoint spherical
systems of generators $(T_1,T_2)$ of  $G$, under the equivalence
relation generated by:
\begin{enumerate}\renewcommand{\theenumi}{\it \roman{enumi}}
\item  Hurwitz equivalence for $T_1$;

\item  Hurwitz equivalence for $T_2$;

\item Simultaneous conjugation for $T_1$ and $T_2$, i.e., for
$\phi \in \Aut(G)$ we let \\ $\bigl(T_1:=(x_{1,1},\dots , x_{r_1,1}),
\ \ \ T_2:=(x_{1,2},\dots, x_{r_2,2})\bigr)$ be equivalent to
\[\bigl(\phi(T_1):=(\phi(x_{1,1}),\dots ,\phi( x_{r_1,1})),
\ \ \phi(T_2):=(\phi(x_{1,2}),\dots, \phi(x_{r_2,2}))\bigr).
\]
\end{enumerate}
Then two surfaces $S$, $S'$ are deformation equivalent if and only
if the corresponding equivalence classes of pairs of spherical
generating systems of $G$ are the same.
\end{theo}

The Hurwitz equivalence is defined precisely in e.g., \cite{P13}. In the cases that we will treat the Hurwitz equivalence is given only by the braid group action on $T_i$ defined as follows. Recall the Artin presentation of the Braid group of $r_1$ strands 
$$
\mathbf{B}_{r_1}:=\langle \gamma_1, \ldots ,\gamma_{r_1-1}| \, \gamma_i \gamma_j = \gamma_j \gamma_i \, \, {\rm for } \, \, |i-j|\geq 2, \, \gamma_{i+1}\gamma_i\gamma_{i+1}=\gamma_i\gamma_{i+1}\gamma_i \rangle.
$$
For  $\gamma_i \in \mathbf{B}_{r_1}$  then:
\[ \gamma_i(T_1)=\gamma_i(v_1, \ldots ,v_{r_1})=(v_1, \ldots ,v_{i+1},v^{-1}_{i+1}v_{i}v_{i+1}, \ldots ,v_{r_1}).
\]

Moreover, notice that, since we deal here with abelian groups only, the braid group action is indeed only by permutation of the elements on the spherical system of generators.


Once we fix a finite abelian group $G$ and a pair of types
$(\tau_1,\tau_2)$ (of size ($r_1$,$r_2$)) of an unmixed ramification
structure for $G$, counting the number of connected components of
$\mathcal{M}_{(G,(\tau_1,\tau_2))}$ is then equivalent to the group
theoretical problem of counting the number of classes of pairs of
spherical systems of generators of $G$ of type $(\tau_1,\tau_2)$
under the equivalence relation given by the action of
$\mathbf{B}_{r_1} \times \mathbf{B}_{r_2} \times \Aut(G)$, given
by:
\begin{equation}\label{eq_count}
    (\gamma_1, \gamma_2, \phi) \cdot (T_1, T_2) := \bigr(\phi(\gamma_1(T_1)),
    \phi(\gamma_2(T_2))\bigl),
\end{equation}
where $\gamma_1 \in \mathbf{B}_{r_1}$, $\gamma_2 \in\mathbf{B}_{r_1}$ and $\phi \in \Aut(G)$, see for more details e.g., \cite{P13}.
\section{Proof of Theorem \ref{thm main}}

Let us consider the group $G:=(\ZZ/2\ZZ)^k$, with $k>>0$ and an integer $l$. We want to give to $G$ many classes of ramification structures of size $(r_1,r_2)=(k(k+1), 2^{l-k+1}+4)$. Since the elements of $G$ have only order two we will produce in the end ramification structure of type $(( 2^{r_1}),( 2^{r_2}))$. 

First let us consider the following elements of $G$
\begin{align*}
v_{1} &= (1,0,\dots,0)     \\
v_{2} &= (1,0,\dots,0)   \\
v_{3} &= (0,1,0,\dots,0)     \\
v_{4} &= (0,1,0,\dots,0)   \\
v_{5} &= (0,1,0,\dots,0)     \\
v_{6} &= (0,1,0,\dots,0)   \\
&\vdots \\
v_{k(k+1)} &= (0,\dots,0,1)  
\end{align*}
and let $T_1:=(v_1, v_2, \ldots ,v_{k(k+1)})$. One can see that $<T_1> \cong G$ and by construction the product of the elements in $T_1$ is $1_G$. Define the set $M:=G\setminus \{0,v_1,\dots, v_{k(k+1)}\} $. We have a bijection
\[ 
M \stackrel{1:1}{\longleftrightarrow} \{n \in \mathbb{N} | n \leq 2^k-k-1 \}=:B.
\]
Call $\varphi \colon B \longrightarrow M$ the bijection map.  Consider $(n_1, \, \ldots ,n_{2^k-k-2})$ a $(2^k-k-1)$-tuple of elements of $B$ whose sum is $2^{l-k}+2$. We define a map 
\[
(n_1, \, \ldots ,n_{2^k-k-1}) \mapsto T_2=\big( \underbrace{\varphi(1), \ldots , \varphi(1)}_{2n_1}, \ldots
,\underbrace{\varphi(2^k-k-1), \ldots , \varphi(2^k-k-1)}_{2n_{2^k-k-1}}\big).
\] 
It holds that $<T_2> \cong G$. Moreover, the product of the elements in $T_2$ is $1_G$, hence $T_2$ is a spherical system of generators for $G$ of size  $2^{l-k+1}$.

By construction $G$ is abelian and all its elements are of order two, therefore the pair $(T_1, \, T_2)$ is a ramification structure for $G$ of the desired type. 
\\

Now we count how many inequivalent ramification structures of this kind we have under the action of the group defined in Theorem \ref{theo_count} and Equation \eqref{eq_count}. First notice that by construction to any tuple $(n_1, \, \ldots n_{2^k-k-2})$ its associated generating vector $T_2$ is in a different braid orbit. Moreover, the choice of $T_1$ implies that any pair $(T_1, T'_2)$ and $(T_1, T''_2)$ are in the same ${\rm Aut}(G)$-orbit if and only if $T'_2=T''_2$.

Hence the number of inequivalent ramification structures is equal to the number of  $(2^k-k-1)$-tuple of positive integers whose sum is $2^{l-k}+2$.  
\begin{quote}
This condition maybe relaxed to the point that only for the elements of 
a basis the entry must be strictly positive and maybe non-negative else.
\end{quote}
This number is known to be
$$
{\frac{r_2}2-1 \choose 2^k-k-2} \quad = \quad {2^{l-k}+1 \choose 2^k-k-2},
$$ 
see e.g., \cite[Section II.5]{F50}.
Let $\nu >0$ be a rational number and let us suppose that $l=(\nu+2) \cdot k$, then using Stirling's approximation of the binomial coefficient - more exactly a corresponding lower bound - we obtain
\begin{equation} \label{eq:bin2}
{2^{l-k}+1 \choose 2^k-k-2} > \frac{(\frac{2^{l-k}+1}{2^k-k-2}-1)^{2^k-k-2}e^{2^k-k-2}}{e\sqrt{(2^k-k-2)}} > (2^{\nu k})^{(2^k-k-2)} \cdot \frac{e^{2^k-k-2}}{e\sqrt{(2^k-k-2)}}> 2^{\nu k(2^k)}.
\end{equation}
Since $e_k = |G| (-2 + \frac12 r_1) ( -2 + \frac12 r_2)$ implies
$2e_k = 2^k \cdot 2^{l-k} \cdot (k^2+k-4)= 2^{(\nu+2)k}(k^2+k-4)$ we have 
$$
(2e_k)^{\frac1{\nu+2}}\cdot \frac{k}{(k^2+k-4)^{\frac1{\nu+2}}}
\quad = \quad k 2^k
$$
Using this,  we obtain for $k$ large enough in the second inequality
\begin{equation} \label{eq:h1}
h > 2^{\nu(2e_k)^{\frac1{\nu+2}}\cdot \frac{k}{(k^2+k-4)^{\frac1{\nu+2}}}}
\quad > \quad 
2^{(e_k^{\frac1{\nu+2}})}
\end{equation}
We can bound further for $k$ large enough
\begin{equation}
2^{(e_k^{\frac1{\nu+2}})}
\quad > \quad
2^{(e_k^{\frac1{2\nu+2}})\frac{\ln e_k}{\ln 2}}
\end{equation}
We use the identity $x^{f(x)}=e^{f(x)\ln x}=2^{f(x)\frac{1}{\ln 2}\ln x}$ to get for all $\alpha<\frac{1}{2}$
\begin{equation*}
h > e_k^{(e_k)^\alpha}
\end{equation*}
if $k$ is large enough, depending on $\alpha$.
This concludes the proof since $e_k$ is proportional to $y_k$.
\hspace*{\fill}$\square$


\end{document}